# Equivalent Effect Function and Fast Intrinsic Mode Decomposition

Louis Yu Lu
E-mail: louisyulu@gmail.com

Abstract: The Equivalent Effect Function (EEF) is defined as having the identical integral values on the control points of the original time series data; the EEF can be obtained from the derivative of the spline function passing through the integral values on the control points. By choosing control points with different criteria, the EEF can be used to find the intrinsic mode function(IMF, fluctuation) and the residue (trend); to fit the curve of the original data function; and to take samples on original data with equivalent effect. As examples of application, results of trend and fluctuation on real stock historical data are calculated on different time scales. A new approach to extend the EEF to 2D intrinsic mode decomposition is introduced to resolve the inter slice non continuity problem, some photo image decomposition examples are presented.

## 1 Introduction

The real time series numeric data from natural phenomena, life science, social and economic systems are mostly nonlinear and non-stationary. The traditional analysis methods like Fourier transform and Wavelet transform presume the linear and stationary of the underlying system, the predefined function bases have little bearing on the physical meaning of the system, and those methods have difficulties to reveal the nature of the real data.

The recent researches tend to use adaptive function bases, Huang et al introduced the IMF as adaptive a posteriori function base in the form of Hilbert spectrum expansion, which has meaningful instantaneous frequency [1][2], The HHT is the Hilbert transform applied to the IMF components. The HHT has broad application to signal processing and times series data analysis in the fields of health, environmental, financial and manufacturing industries [4][8].

After a brief introduction to HHT and IMF, the problems with HHT and related methods will be analyzed. Then EEF definition and the steps of finding the EEF will be presented. Some real world applications on 1D (stock historical data trend and fluctuation), 2D (photo image decomposition)  will be showcased.

The conclusion and future study direction will be discussed by the end.

## 2  Intrinsic Mode Function (IMF)

Huang et al analyze the requirement of meaningful instantaneous frequency on Hilbert transformation, and introduced the adaptive function base as IMF [1][2], which satisfies the two conditions:
1) In the whole data set, the number of extrema and number of zero crossings must equal or differ at most by one;
2) At any point, the mean value of the envelope defined by the local maxima and the envelope defined by the local minima is zero.

The first condition is apparently necessary for oscillation data; the second condition requires the symmetric upper and lower envelopes of IMF, as the IMF component is decomposed from the original data; it is quite challenging to find the real envelopes because of nonlinear and non-stationary nature in the data. Only a few functions have the known envelopes, for example, the constant amplitude sinusoidal function.

All the IMF has the meaningful instantaneous frequency defined by the Hilbert transformation [3]. The different IMF components of a time series data reflects the oscillation of the data on different time scales, while the residue components reflect the data trend.

The problem with the IMF definition is no guaranty of zero integral value on the whole data range of IMF, it is less likely the summation be zero on local interval. In a physical sense related to electricity, the IMF still contains direct current components.

## 3  Empirical Mode decomposition (EMD)

To break up the original data into a series of IMF, Huang et al invented the EMD method [1][4], the idea is to separate the data into a slow varying local mean part and fast varying symmetric oscillation part, the oscillation part becomes the IMF and local mean the residue, the residue serves as the input data again for further decomposition, the process repeats until no more oscillation can be separated from the residue of higher frequency mode.

On each step of the decomposition, since the upper and lower envelope of the IMF is unknown initially, a repetitive sifting process is applied to approximate the envelopes with cubic spline functions passing through the extrema of the IMF, the data function serves as the initial version of IMF, and the refining IMF is calculated as the difference between the previous IMF version and mean of the envelopes, the process repeats until the predefined stop condition satisfied. The residue is then the difference between the data and the improved IMF.

There are four problems with this decomposition method:
1) The spline (cubic) connecting extrema is not the real envelope, some unwanted overshoot may be introduced by the spline interpolation, the resulting IMF function does not strictly guarantee the symmetric envelopes [1]. Higher order spline does not in theory resolve the problem.
2) The different stopping condition value results different set of IMF[1], picking up the right value is a more subjective matter, making the results unpredictable.
3) The repetitive sifting process is time consuming, because the problem in item 2, more sifting does not guaranty better results. Some researches try to improve the performance, the study in reference only speeds up the envelopes building process; problem 1 and 2 still exist.

4) Some riding wave on steep edge of IMF are missing [1], the resulting IMF may not accurately reflect the true physical nature of the data.

# 4   Previous Study of Fast Intrinsic Mode Decomposition Methods by the Author

Previously, the author had developed two methods to improve the predictability and speed up the algorithm of decomposing the intrinsic modes.

The first is the sawtooth transform method[5], the approach is to transforms the original data function into a piecewise linear sawtooth function(or triangle wave function), then directly constructs the upper envelope by connecting the maxima and construct lower envelope by connecting minima with straight line segments in the sawtooth space, the IMF is calculated as the difference between the sawtooth function and the mean of the upper and lower envelopes. The results found in the sawtooth space are reversely transformed into the original data space as the required IMF and envelopes mean. This decomposition method processes the data in one pass to obtain a unique IMF component.

The problem with the method is the leaking of high frequency ripples into the residue function.

The second approach[6] is to iteratively adjust the control points on the data function corresponding to the extrema of the refining IMF, the control points of the residue function are calculated as the median of the straight line segments passing through the data control points, the residue function is then constructed as the cubic spline function of the median points. The initial residue function is simply constructed as the straight line segments passing through the inflexion points of the data function. The refining IMF is the difference between the data function and the improved residue function. The IMF found reveals all the riding waves in the whole data set.

This result of the method still suffers the problem with the definition of IMF, the IMF components contain direct current and the residues do not reflect the true trends of the data.

# 5   Definition of Equivalent Effect Function (EEF)

Let $f(t)$ be the original time series data function and $F(t)$ the integral, $f(t)$ has $m$ control points:

$$F(t) = \int_{t_0}^{t} f(x)\,dx \qquad t_0 \leqslant t < t_m$$

Let $e(t)$ be a function on the same domain as $f(t)$, and $E(t)$ the integral:

$$E(t) = \int_{t_0}^{t} e(x)\,dx \qquad t_0 \leqslant t < t_m$$

The $e(t)$ is defined as the EEF of $f(t)$ satisfying the conditions:

$$E(t_i) = F(t_i) \qquad 0 \leqslant i < m$$

The conditions force $e(t)$ and $f(t)$ having the accumulated equivalent effect on the control points. The difference between $f(t)$ and $e(t)$ is defined as:

$$d(t) = f(t) - e(t) \qquad t_0 \leqslant t < t_m$$

The difference function $d(t)$ has the property of zero integral value on control points and zero summation on the whole data range.

# 6 Finding EEF and its Representation with Spline

With the definition in chapter 5, it is straight forward to obtain the EEF with the following steps:
1) Calculate the integral of *f(x)* on the control points $F(t_i)$;
2) Find a spline function *E(t)* passing through the control points $(t_i, F(t_i))$;
3) The EEF is calculated as the derivative of *E(t)*:
$$e(t) = E'(t) \qquad t_0 \leq t < t_m$$

By choosing different order of spline function to represent *E(t)*, *e(t)* has different smoothness. The existence of proven method of finding spline function on different order has effect on the cost and reliability of the algorithm.
1) Using piecewise linear function for *E(t)*, *e(t)* becomes non continuous step function, it is not quite useful in the context of this study.
2) Using quadratic spline function for *E(t)*, *e(t)* becomes piecewise linear function. As the even order spline function oscillate around the control points, a method is introduced in chapter 7 to solve the problem.
3) Using cubic spline function for *E(t)*, *e(t)* becomes quadratic spline function, the second derivative becomes non continuous. The algorithm to find the cubic spline is very popular and efficient.
4) Using quartic spline function for *E(t)*, *e(t)* becomes cubic spline function, the second derivative becomes continuous. As the even order spline function oscillate around the control points, a method is introduced in chapter 7 to solve the problem.
5) Using quintic spline function for *E(t)*, *e(t)* becomes quartic spline function, the third derivative becomes continuous. The quintic spline algorithm developed J.G. Herriot and C. H. Reinsch[7] is very efficient, this algorithm also allows specifying first and second derivatives on any knot besides the boundaries. So far quintic spline function for *E(t)* has best results for the experimenting of intrinsic mode decomposition.

# 7 Even Order Spline Function

Finding the spline is in general turned into solving a system of linear equation, the number of unknown must be equal to the number of independent equations.
Given *n* sample points:
$$(x_i, y_i) \qquad 0 \leq i < n$$
The even order (2 or 4) spline function passing through the samples needs to be found. The traditional way links segment of polynomial functions on the sample points, but the number of boundary constraints equations are not even number, the non-symmetric constraints on the two ends may cause the oscillation of the spline function.
To bring a remedy to the problem, the segment functions are rearranged passing through the samples and connecting knots in between sample points, resulting *n-1* knots :
$$(u_j, v_j) \qquad 0 \leq j < n-1$$
Where:
$$u_j = \alpha x_j + (1-\alpha) x_{j+1} \qquad 0 < \alpha < 1$$
This value of *α* can simply be chosen as 0.5 to position the knot in the middle of two samples, With $x_0$ connecting $u_0$, $u_j$ connecting $u_{j+1}$, and $u_{n-2}$ connecting $x_{n-1}$, there are total *n* segments, that requires *n* polynomial function to build the spline.

1) For quadratic spline, each segment function has 3 coefficients as the unknown variables:
$$y = a + bx + cx^2$$
The spline contains **3n** unknowns, the value $y_i$ on the samples provide **n** equation, the equality of value and first derivative from the two neighbouring functions on the knots provides **2(n-1)** equations, 2 symmetric boundary constraints(specific value on first derivative or second derivative on each end) provide **2** extra equations, the total numbers equations:
$$n + 2(n-1) + 2 = 3n$$
This is the exact number of equations to resolve the unknown variable in the quadratic spline function.

2) For quartic spline, each segment function has 5 coefficients as the unknown variables:
$$y = a + bx + cx^2 + dx^3 + ex^4$$
The spline contains **5n** unknowns, the value $y_i$ on the samples provide **n** equation, the equality of value and up to third derivative from the two neighbouring functions on the knots provides **4(n-1)** equations, 4 symmetric boundary constraints(specific value on first derivative and second derivative on each end, or not a knot condition)provide **4** extra equations, the total numbers equations:
$$n + 4(n-1) + 4 = 5n$$
This is the exact number of equations to resolve the unknown variable in the quartic spline function.

# 8  Fast Intrinsic Mode Decomposition(FastIMD) Method with EEF

The EEF generates smooth function passing through the control points. With the result of previous study[6], the polarity of the riding wave is separated by the inflexion point on the data function, choosing the control points in between the inflexion points will generate an EEF which smooth out the original data. The EEF will serve as the trend or residue function, the difference function serves as the IMF or fluctuation function.

The following steps describe the algorithm of the Fast Intrinsic Mode Decomposition
1) Calculate the first derivative of the data function, find the extrema on the derivatives which are the inflexion points;
2) Connect the inflexion points with straight line segments, extend the first segment to the first sample and last segment to the last sample, this piecewise linear function is the estimate trend function;
3) Calculate the difference between the original data and the estimate trend function obtained in step 2, find the extrema on the difference, these extrema will serve as the control points of the EEF;
4) Calculate the integral values on the control points obtained in step 3, build spline function on these integral values;
5) Calculate the derivative of the integral spline function of step 4, the resulting EEF is the trend(residue) function;
6) Calculate the difference between the data function and the trend function, the result is the fluctuation function or IMF;

7) Take trend function as the data function and repeat from step 1 to find another trend and fluctuation component at lower frequency, the process stops when the extrema count goes down to 2.

This algorithm is fast and reliable as there is no sifting or complicated calculation involved. Regarding the selection of control points, as long as they are located between the inflexion points, the oscillation frequency on the successive trend will be reduced, for example the control point can be simply chosen as the middle points between the neighbouring inflexion points. But the experiments generate less number of IMF components with method described in step 3 than alternative ones.

Another way to choose the control points is to use the extrema on the data function, but it will miss the riding waves on steep edge of IMF.

## 9 Boundary Constraints

With quintic spline[7] on the integral data, the specification of boundary constraints can be quite flexible. On either end, there can be no constraints, or specified first derivative only, or specified both first and second derivatives.

The fist derivative corresponds to the value on EEF, and second derivative to the first derivative on EEF, these allow to chain piecewise EEF together with continuity on first derivative, which can benefit the processing of long duration time series data.

## 10 Curve Fitting with EEF

EEF can also serves the purpose of curve fitting, by choosing the the control points as the union of extrema and inflexion points of the original data, the resulting EEF fits the original data very nicely, and can be used as data compression.

The method is to calculate the first derivative of the data function, convert the result into absolute values, then find the extrema on the absolute values. The extrema on absolute values of the first derivative are the inflexion points and extrema on the original data.

With the control points in hand, the fitting spline is obtained similar to step 4 and 5 in chapter 8.

With FastIMD algorithm, all the fluctuation and trend function are polynomial except the first fluctuation function with the highest frequency, it can be fitted to become polynomial function.

## 11 Data Sampling with EEF

Data sampling can be done by choosing a subset of the sample points either periodically or selectively as the control points, then build the EEF with the method similar to step 4 and 5 in chapter 8.

The data sampling with EEF has the advantage of equivalent effect comparing to the traditional averaging method. It can be used for data compression or kind of low pass filtering.

## 12 Trend Analysis of Historical Stock Transactions

The FastIMD algorithm described in chapter 8 is used to develop a financial data analysis tool 'EEF Trend' coded with ActionScript and posted on BlackBerry App World.

'EEF Trend' allows the user to find the stock symbol of a company, obtain the current stock quote, retrieve the historical daily trading data of the company from the Internet (Yahoo Finance) and graphically display the trend and fluctuation with varying time scales from the fast moving short term trend to the slow moving long term general trend.
In the illustrated graphs, the legend corresponds the following terms:
Original -- the original data associated with the time;
Trend n -- the trend data of mode n, the bigger the number n, the less fluctuation contained, the last trend with the biggest number n contains no fluctuation;
Difference n -- the difference between the original data and the trend n;
Fluctuation n -- the difference between the previous trend and trend n.
The data represent the daily stock closing value of 'Nasdaq Composite' for the past five years, total 11 modes have been generated. Mode 1, 2, 3, 5, 8, 11are displayed bellow.

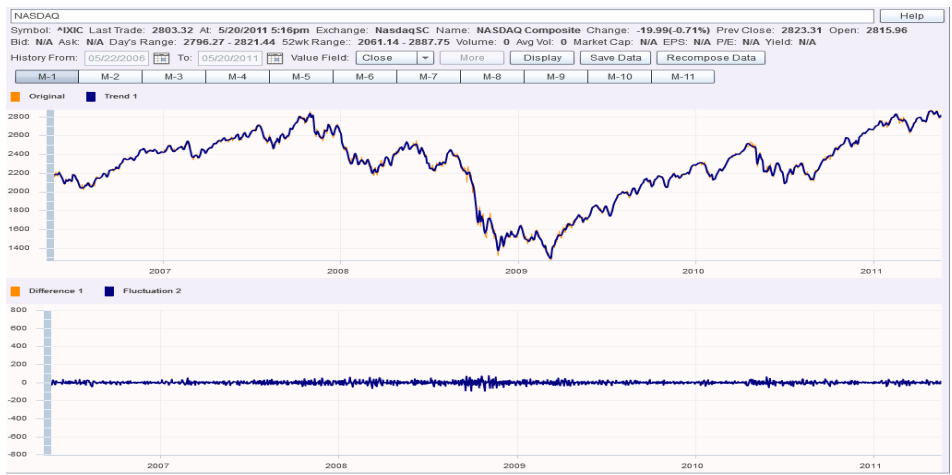

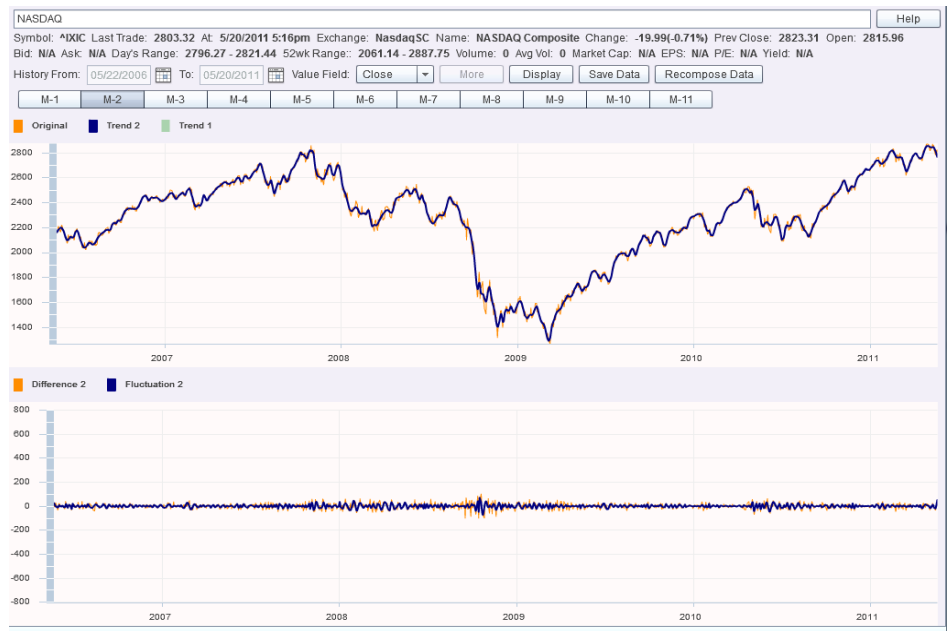

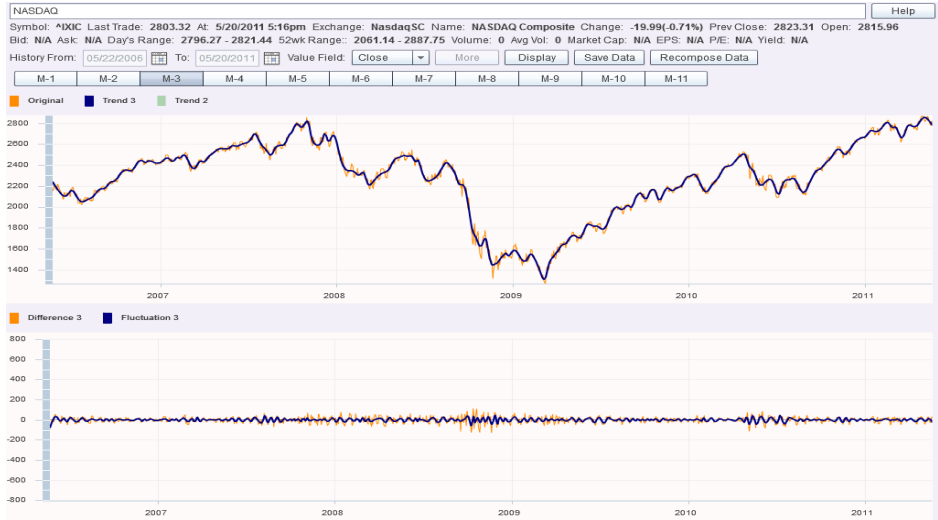

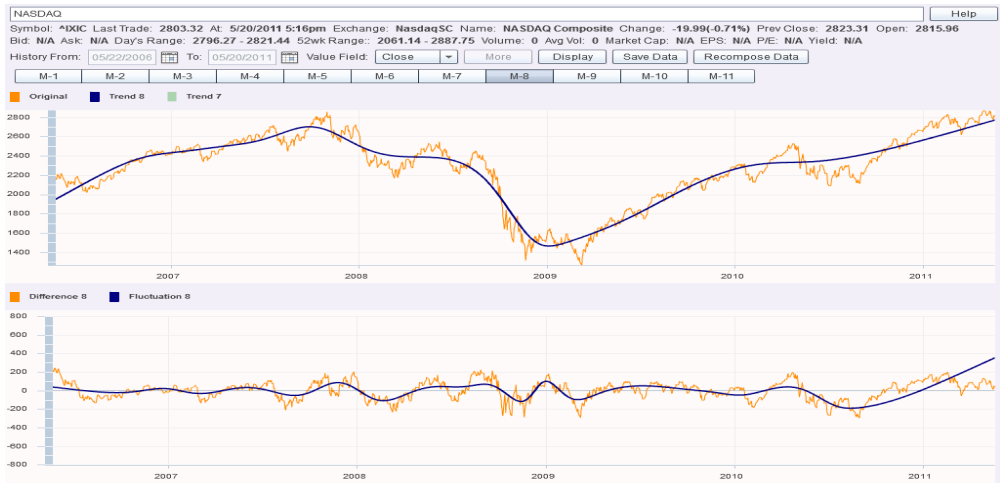

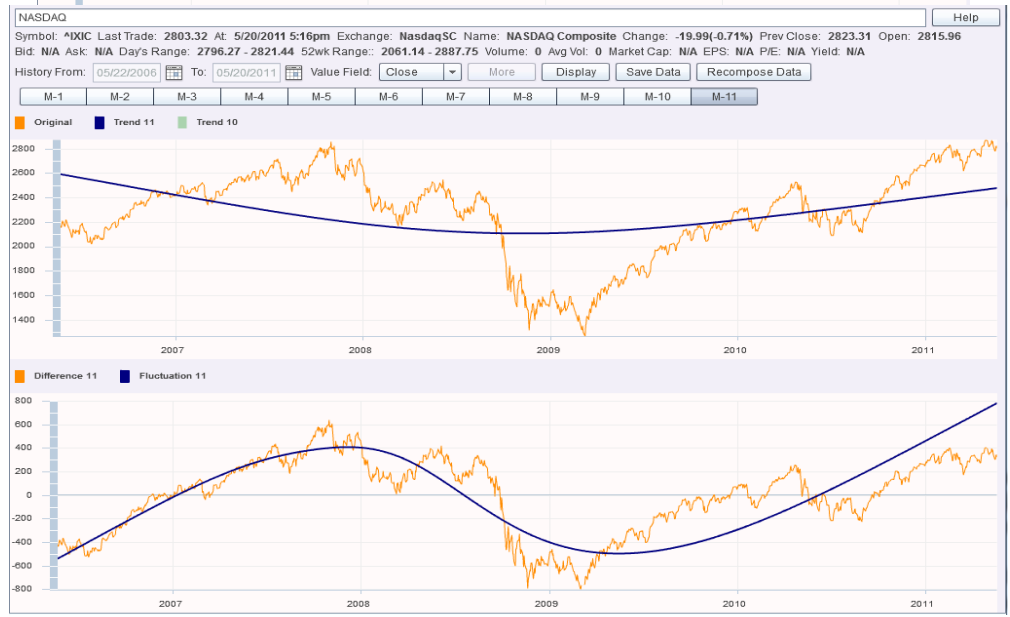

## 13  FastIMD 2D Extension and Image Processing

The FastIMD algorithm is extended to 2D for image data by treating the data as a collection of 1D slices, the slice can be horizontal scan line or vertical scan line. This similar approach was taken by the pseudo-two-dimensional EMD[8], but the pseudo-two-dimensional EMD has shortcomings, and one of them is the inter slice discontinuity.

To solve the inter slice discontinuity problem, a slight modification of the FastIMD is applied. Since the contiguous slices are continuous on the data function, so the contiguous slice trend functions should be continuous as well. Instead of using the piecewise linear function as the estimate trend function as in step 3 of chapter 8, the trend function of the previous slice is used as the estimate trend function for the current slice, the other steps is similar to the 1D FastIMD in chapter 8.

The following steps describe the FastIMD 2D Extension:
1) For the first slice or scan line, use the regular 1D FastIMD described in chapter 8 to find the trend functions;
2) For the rest of slices, use the trend function of the previous slice as the estimate trend function, following steps 3, 4, and 5 in chapter 8 to calculate the trend function for the current slice;
3) Assemble of slice trend functions to form the 2D trend function in one scan direction;
4) Take the 2D trend function as data, start the slice in perpendicular direction and go to step 1 to find trend function in another direction, calculate the difference between the data function and the trend function with two directions scan to construct the fluctuation function.
5) Repeat the steps above to find the trend and fluctuation 2D function with lower frequency.

Two pictures are used to show case the algorithm, the first is the well known Lena portrait for image processing, the second is a regular picture of the author and his son. The colour image is first converted into black and white grey scales, then the grey scales image is applied with the decomposition algorithm, first scan in horizontal direction, then in vertical direction.

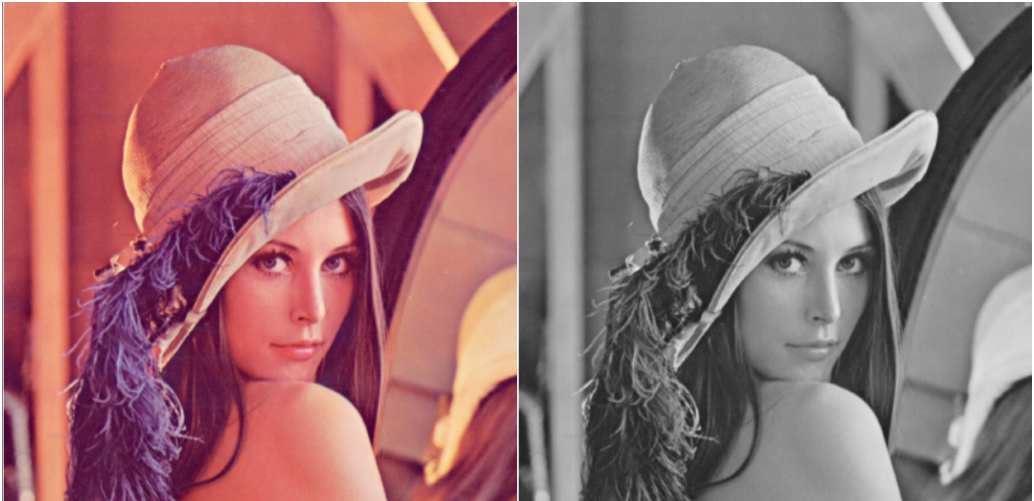

The standard Lena portrait for image processing

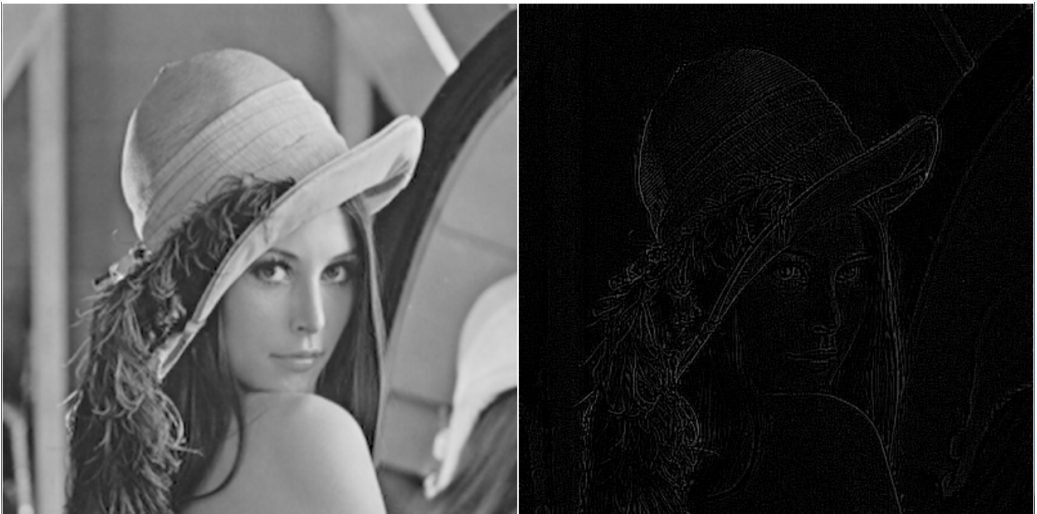

Mode-1, left: 2D trend; right: 2D fluctuation

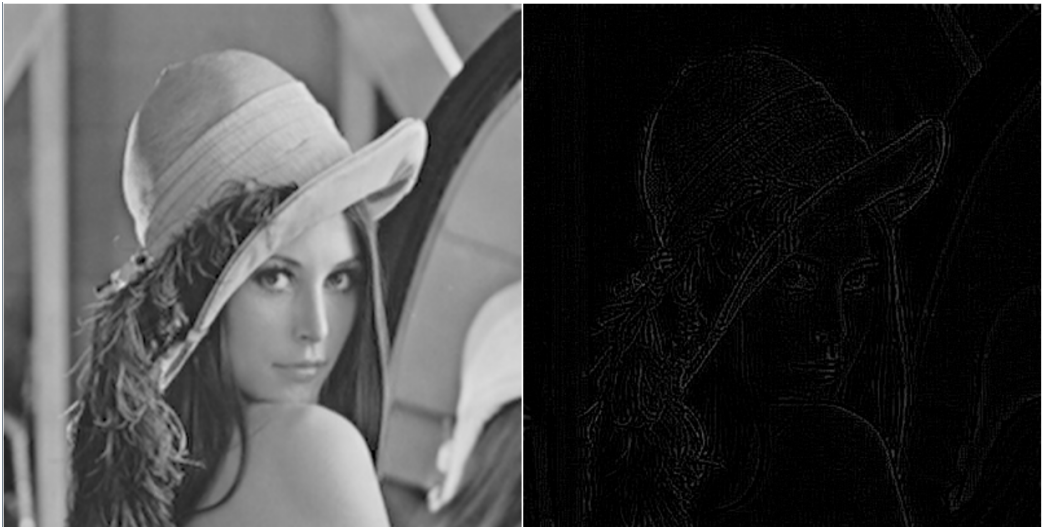

Mode-2, left: 2D trend; right: 2D fluctuation

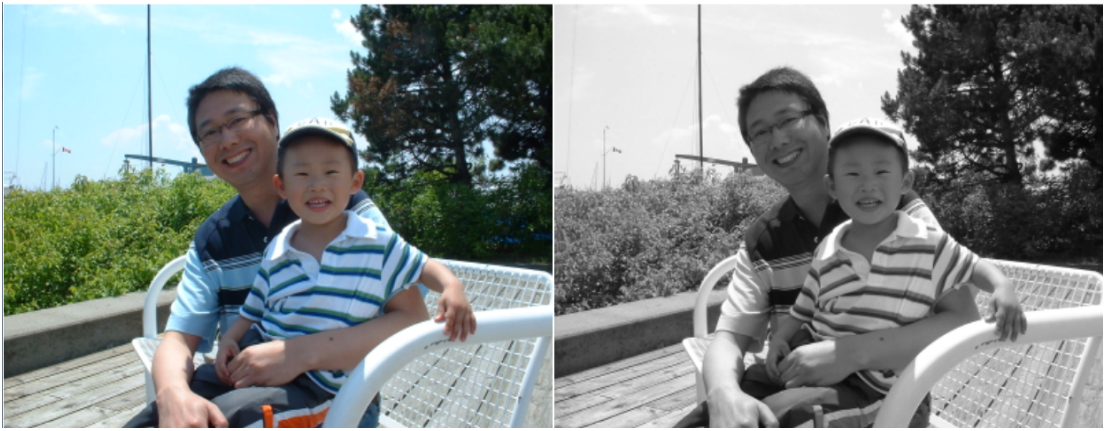
Louis and his son Tylor (taken in June 2008)

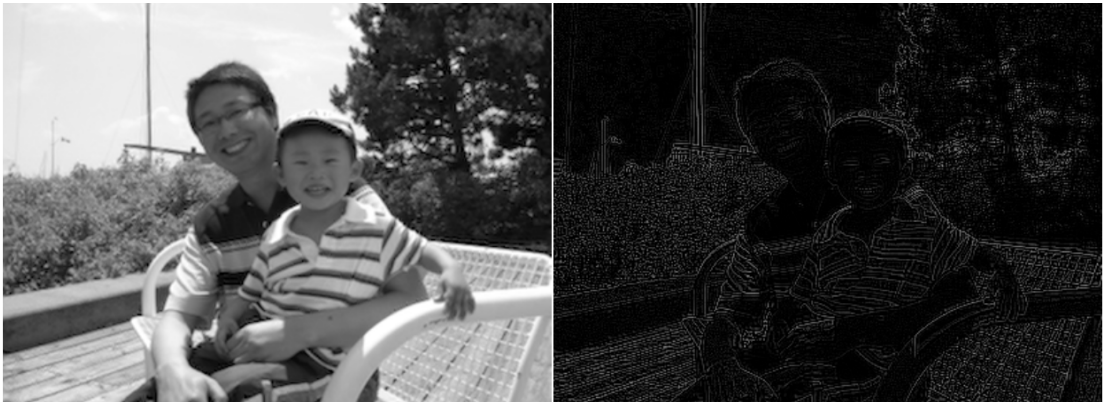
Mode-1, left: 2D trend; right: 2D fluctuation

## 14  Conclusion

The EEF introduced in this paper serves as a fast method to find the IMF, it can also be used for curving fitting and data sampling. The experimenting results for 1D data and 2D image data processing are quite promising.

The more rigorous theoretical background of EEF than EMD and better resulting IMF will position it as serious analysis method complementary to traditional Fourier analysis and Wavelet analysis methods.

The FastIMD method presented is efficient and generating satisfactory results in a practical sense, the theoretical criteria for the best results still need to be studied and clarified.

The proposed method of resolving quartic spline is base on a general method of solving the linear systems equations, a more efficient and reliable method comparable to the quintic spline[7] is desirable. Besides polynomial spline function for EEF, other interpolation functions are worth exploring.

In addition to the 2D EEF extension with 1D slice, a direct 2D EEF method are expected to be developed, it is still not clear whether it can generate satisfactory result by building a surface function having same 2D integral value on the extrema and saddle points of the original data.